\documentclass[11pt]{amsart}
\usepackage{amsmath, amssymb, amsthm}
\usepackage{bm}
\usepackage[numbers]{natbib}
\usepackage{textcomp}
\usepackage{bibentry}
\usepackage{etoolbox}
\usepackage{hyperref}
\hypersetup{
  pdfauthor={Tianyuan Workshop},
  pdftitle={Open Problems in Mathematical Logic},
  pdfsubject={Open Problems},
  colorlinks=true,
  linkcolor=blue,
  citecolor=blue,
  urlcolor=blue
}
\usepackage{lmodern}
\usepackage{tikz}
\usepackage{relsize}
\usepackage{parskip}
\usepackage{enumerate}
\usepackage{textcomp} 
\theoremstyle{definition}
\newtheorem{theorem}{Theorem}[section]
\newtheorem{problem}{Problem}[section]

\newtheorem{proposition}[theorem]{Proposition}
\newtheorem{example}[theorem]{Example}
\newtheorem{remark}[theorem]{Remark}
\newtheorem{definition}[theorem]{Definition}
\newtheorem{fact}[theorem]{Fact}
\newtheorem{notation}[theorem]{Notation}
\def\graph{\mathrm{GRAPH}}
\def\itx#1{{\mbox{\textrm{#1}}}}

\newcommand{\abstractauthor}[1]{\noindent\textbf{#1}\par\vspace{0.3\baselineskip}}
\newcommand{\abstractbody}[1]{\noindent #1 \par \vspace{\baselineskip}}

\makeatletter
\renewcommand{\bibsection}{\vspace{\baselineskip}\noindent\textbf{References}}
\makeatother

\begin{document}

\title[Open Problems in Mathematical Logic]{Open Problems in Mathematical Logic}
\author[2026 Tianyuan Workshop in Definability and Computation]{George Barmpalias}  
\author[]{Su Gao}  
\author[]{Jialiang He}  
\author[]{Takayuki Kihara}  
\author[]{Andre Nies}  
\author[]{Theodore Slaman}  
\author[]{Chieu-Minh Tran}  
\author[]{Daniel Turetsky}  
\author[]{Philip Welch}  
\author[]{Liang Yu}  
\author[]{Hang Zhang}


\begin{abstract}
These open problems were presented in the Problem Sessions held during the Tianyuan Workshop on Definability and Computation, June 22-26, 2026. The problems are organized into sections named after their contributors, in the order of their presentations during the workshop. Notes were taken and compiled by Wei Dai, Xiangxi Hu, Yingying Jiang,  Ruiwen Li, Tianhao Wang, Xu Wang, and Jie Zou.
\end{abstract}
\maketitle
\tableofcontents

\section{George Barmpalias}

The following problems concern Turing ideals generated by sufficiently
random reals.

\begin{problem}[Levin, private communications]
        Are the following true?
    \begin{enumerate}
        \item[(1)] Can every chain of random reals in the Turing degrees
        be combined into a single random real?

        \item[(2)] Can every Turing ideal generated by random reals be
        generated by the columns of a single random real?
    \end{enumerate}
\end{problem}

\begin{notation}
If $\alpha\in 2^\omega$, let $\alpha^{(i)}$ denote the real obtained from
$\alpha$ by deleting the digits $\alpha(k)$ at positions $k$ that are
divisible by $2^i$.
\end{notation}

\begin{itemize}
    \item A \emph{chain} is a sequence $(x_i)_{i\in\mathbb{N}}$ of
    reals such that
    \[
        x_i<_T x_{i+1}
    \]
    for every $i\in\mathbb{N}$.

    \item If $B$ is a chain of reals, the Turing ideal generated by
    $B$ is
    \[
        [B]:=\{z:\exists x\in B,\ z\leq_T x\}.
    \]
    In particular, for a chain $(x_i)_{i\in\mathbb{N}}$,
    \[
        [(x_i)]
        :=
        \{z:\exists i,\ z\leq_T x_i\}.
    \]

    \item For each real $\alpha$, define
    \[
        [\alpha]
        :=
        \{z:\exists i,\ z\leq_T\alpha^{(i)}\}.
    \]

    \item The notation $z\approx x$ means that $z$ and $x$ differ
    on finitely many bits.

    \item The real
        $\bigoplus_i x_i$
    is the infinite join of all the reals $x_i$.
\end{itemize}

\medskip
\begin{proposition}[Barmpalias--Zhang, unpublished]
\textit{If $(x_i)_{i\in\mathbb{N}}$ is a sequence of reals and $\bigoplus_{i\leq n} x_i$ is random for each $n$, then there are reals
$z_i\approx x_i$ so that $\bigoplus_i z_i$ is random. The same conclusion holds for $k$-randomness and weak $k$-randomness for each $k\in \mathbb{N}$.}
\end{proposition}

This proposition gives a partial positive answer to Question~1(1) under
the additional assumption of relative randomness. Question~1(2) asks for
the stronger conclusion that the ideal generated by a random chain is
exactly of the form $[\alpha]$ for some random real $\alpha$.

\begin{problem}[Levin]
    Is it true that, for every $X\subseteq 2^\omega$ with $\mu(X)=1$,
    there exists $Y\subseteq 2^\omega$ with $\mu(Y)=1$ such that, for
    every chain $B\subseteq Y$, there is some $\alpha\in X$ satisfying
    $[B]=[\alpha]$?
  \end{problem}

The class $\mathrm{ML}(\emptyset')$
of $\emptyset'$-random reals, equivalently the class of $2$-random reals,
is a special case.

\begin{problem}
    Is there a chain
    $B\subseteq\mathrm{ML}(\emptyset')$
    such that
    $[B]\neq[\alpha]$
    for every $\alpha\in\mathrm{ML}(\emptyset')$?
\end{problem}

There are two obstacles to constructing such a chain. If $\alpha$ is
$2$-random, then the ideal $[\alpha]$ has the following properties:
\begin{enumerate}
    \item[(i)] For every $x\in[\alpha]$, the ideal $[\alpha]$ contains an
    $x'$-random real.

    \item[(ii)] The jumps of the columns $\alpha^{(i)}$ form a chain and,
    in particular, are pairwise distinct.
\end{enumerate}

\begin{problem}
    Toward a negative answer to Levin's Problem 1.2, are the following true?
    \begin{enumerate}
        \item[(1)] Are there $2$-random reals $x<_T z$ such that
        \[
            x'\equiv_T z'?
        \]

        \item[(2)] Are there $2$-random reals $x<_T z$ such that
        \[
            \left|
                K(x\upharpoonright n)-K(z\upharpoonright n)
            \right|
            =
            O(1)?
        \]

        \item[(3)] Is there a Turing ideal $B$ generated by a chain of
        $2$-random reals such that $B$ contains no $x$-random real for
        any $x\in B$?

        \item[(4)] What are the jumps of $2$-random reals?
        See \cite{BDN}. 
    \end{enumerate}
\end{problem}

\begin{problem}
    Is there an injective one-way function on the reals? More precisely,
    is there a partial computable real function $f$ such that
    \begin{enumerate}
        \item[(i)] $f$ is injective and preserves randomness, 

        \item[(ii)] the domain of $f$ has positive Lebesgue measure, and 

        \item[(iii)]
        $f(x)<_T x$
        for a positive-measure set of reals $x$ in the domain of $f$?
    \end{enumerate}
\end{problem}

Such a function can be used to produce a chain of $2$-random reals with
the same jump in the Turing degrees.

\vskip -10pt
\begin{bibentry}{9}

\end{bibentry}

\bigskip
\bigskip

\section{Andre Nies}
\subsection*{Oligomorphic groups}

Let $M$ be a countably infinite $\omega$-categorical structure in a
finite first-order signature, and let
\[
    G=\operatorname{Aut}(M),
\]
equipped with the topology of pointwise convergence. Let
$\operatorname{Aut}(G)$ denote the group of continuous automorphisms
of the topological group $G$, and let
\[
    \operatorname{Out}(G)
    :=
    \operatorname{Aut}(G)/\operatorname{Inn}(G)
\]
be its outer automorphism group.

\begin{problem}
    Let $M$ be an $\omega$-categorical structure in a first order finite language. Is the group $\operatorname{Out}(G)$  finite?
\end{problem}

For background see \cite{NiesPaolini}, where it is shown that this group is  countable. For some classical examples like $M=(\mathbb{Q},<)$ \cite{Rubin}, the statement of the problem is true; in this case ${\rm Out} (G)$ has two elements, given by the identity and conjugation by $x \mapsto -x$.    

\subsection*{Computability Theory}

Assume throughout that $A,B\subset\mathbb{N}$. One writes 

\[
A \le_{SJT} B
\]

if for every   order function $h$ there exists a uniformly
$B$-c.e. trace $(T_n)$ satisfying that

\[
|T_n|\le h(n)
\]

and

\[
J^A(n)\downarrow
\Longrightarrow
J^A(n)\in T_n.
\]
Also recall that 
\[
A\le_{LR} B
\iff
{\rm MLR}^B\subseteq{\rm MLR}^A.
\]

\begin{problem}
	Does $\le_{SJT}$ imply $\le_{LR}$?
\end{problem}

 We note that  $A\le_{SJT}\emptyset$ implies $A \le_{LR}\emptyset$, namely, every strongly jump traceable set is $K$-trivial, equivalently, low for ML-random.  For background and references, see \cite{GNT}. The authors there give several   characterisations of $\le_{SJT}$ on the $K$-trivials; for instance,  for $A,B$ $K$-trivial,

\[
A\le_{SJT} B
\iff
\forall\,Y\, \omega\text{-c.a. and }Y\in{\rm MLR}
\Bigl(
A\le_T B\oplus Y
\Bigr).
\]

Let Rec denote the class of recursive sets.

\begin{problem}
    Is there a $\Pi^0_1$ class $C$ such that $|C\cap {\rm Rec}|=\infty$ and $C\cap {\rm Rec}$ is maximally almost disjoint among the recursive sets?
\end{problem}

\vskip -12pt

\begin{bibentry}{9}

\end{bibentry}

\bigskip
\bigskip

\section{Takayuki Kihara}
\abstractbody{

The following problems concern a realizability-theoretic refinement of
many-one and Wadge reducibility.

For sets $A,B\subseteq\omega^\omega$, define
\[
    A\leq_m B
\]
if there is a computable function
\[
    \theta\colon\omega^\omega\to\omega^\omega
\]
such that
\[
    x\in A
    \quad\Longleftrightarrow\quad
    \theta(x)\in B
\]
for every $x\in\omega^\omega$. Replacing ``computable'' by
``continuous'' gives Wadge reducibility, denoted by
\[
    A\leq_W B.
\]

For formulas, the truth-preserving map on instances is not sufficient.
The witnesses for the two formulas must also be uniformly convertible.

Let $\alpha\Vdash\varphi(x)$ mean that $\alpha$ is a witness, or
realizer, for the truth of $\varphi(x)$. The witness relation is defined
inductively from the logical structure of $\varphi$; see
\cite[Definition~2.10]{KiharaHierarchy}.

A formula $\varphi$ is realizability-theoretically many-one reducible
to a formula $\psi$, written
\[
    \varphi\leq_m\psi,
\]
if there are computable functions $\eta,r^{-},r^{+}$ such that
\begin{enumerate}
    \item[(1)]
    \[
        \varphi(x)
        \quad\Longleftrightarrow\quad
        \psi(\eta(x));
    \]

    \item[(2)] If $\alpha\Vdash\varphi(x)$, then
    \[
        r^{-}(\alpha,x)\Vdash\psi(\eta(x));
    \]

    \item[(3)] If $\beta\Vdash\psi(\eta(x))$, then
    \[
        r^{+}(\beta,x)\Vdash\varphi(x).
    \]
\end{enumerate}
This is also called \emph{Levin reducibility}; see
\cite[Definition~2.13]{KiharaHierarchy}. Its continuous analogue gives
the corresponding realizability-theoretic version of Wadge
reducibility.

Existentially quantified variables may depend on all preceding
universally quantified variables. For example, a witness for
\[
    \forall a\,\exists b\,\forall c\,\exists d\,
    \forall e\,\exists f\,\forall g\,
    \varphi(a,b,c,d,e,f,g,x)
\]
consists of three functions
\[
    F_0\colon\omega\to\omega,\qquad
    F_1\colon\omega^2\to\omega,\qquad
    F_2\colon\omega^3\to\omega
\]
such that
\[
    \forall a\,\forall c\,\forall e\,\forall g\,
    \varphi\bigl(
        a,F_0(a),
        c,F_1(a,c),
        e,F_2(a,c,e),
        g,x
    \bigr).
\]
Thus, $F_0(a)$ witnesses $b$, $F_1(a,c)$ witnesses $d$, and
$F_2(a,c,e)$ witnesses $f$.

In addition to the ordinary quantifiers $\exists$ and $\forall$, define
\[
    \exists^\infty n\,\varphi(n)
    \quad\Longleftrightarrow\quad
    \forall m\,\exists n\geq m\,\varphi(n),
\]
and
\[
    \forall^\infty n\,\varphi(n)
    \quad\Longleftrightarrow\quad
    \exists m\,\forall n\geq m\,\varphi(n).
\]

A finite sequence
\[
    \overline{Q}
    \in
    \{\exists,\forall,\exists^\infty,\forall^\infty\}^{<\omega}
\]
is called a \emph{quantifier pattern}. If
\[
    \overline{Q}=Q_0Q_1\cdots Q_\ell,
\]
then a formula of the form
\[
    Q_0n_0\,Q_1n_1\cdots Q_\ell n_\ell\,
    \theta(n_0,\ldots,n_\ell,x),
\]
where $\theta$ is bounded, is called an
$\overline{Q}$-formula.

For each quantifier pattern $\overline{Q}$, let
\[
    \langle\overline{Q}\rangle
\]
denote a fixed $\overline{Q}$-complete formula.

Consider
\[
    \mathrm{Bdd}
    :=
    \left\{
        x\in\omega^\omega:
        \exists n\,\forall k\,
        x(k)\leq n
    \right\},
\]
and
\[
    \mathrm{FIN}
    :=
    \left\{
        x\in\omega^\omega:
        \exists n\,\forall m\geq n\,
        x(m)=0
    \right\}.
\]

Although both are classically $\Sigma^0_2$ properties, their witness
structures are different. The first has the quantifier pattern
$\exists\forall$, whereas the second has the pattern
$\forall^\infty$.

The realizability-theoretic classification of
$\Sigma^0_2$ quantifier patterns consists of exactly three degrees:
\[
    \langle\forall^\infty\rangle
    <_m
    \langle\forall^\infty\forall\rangle
    <_m
    \langle\exists\forall\rangle.
\]
The abstract treatment of this classification is given in
\cite{KiharaManyOne}.

The more concrete analysis of the $\Sigma^0_3$ and $\Pi^0_3$ levels is
given in \cite{KiharaHierarchy}. The known numbers of equivalence
classes are as follows:
\begin{enumerate}
    \item[(1)] There is exactly one many-one equivalence class of $\Pi^0_2$ quantifier patterns,
    represented by
    \[
        \forall\exists.
    \]

    \item[(2)] There are exactly three many-one equivalence classes of $\Sigma^0_2$ quantifier
    patterns, represented by
    \[
        \exists\forall,\qquad
        \forall^\infty\forall,\qquad
        \forall^\infty.
    \]

    \item[(3)] There are exactly three many-one equivalence classes of
    $\Sigma^0_3$ quantifier patterns, represented by
    \[
        \exists\forall\exists,\qquad
        \forall^\infty\exists^\infty,\qquad
        \forall^\infty\exists.
    \]

    \item[(4)] There are exactly five many-one equivalence classes of
    $\Pi^0_3$ quantifier patterns, represented by
    \[
        \forall\exists\forall,\qquad
        \exists^\infty\forall^\infty\forall,\qquad
        \exists^\infty\forall,\qquad
        \forall\forall^\infty\forall,\qquad
        \forall\forall^\infty.
    \]
\end{enumerate}

The corresponding classifications at the fourth level are still in
progress.

\begin{problem}
    Determine the equivalence classes of quantifier patterns at the
    fourth level of the arithmetical hierarchy:
    \begin{enumerate}
        \item[(1)] How many $\equiv_m$-equivalence classes of
        $\Sigma^0_4$ quantifier patterns are there?

        \item[(2)] How many $\equiv_m$-equivalence classes of
        $\Pi^0_4$ quantifier patterns are there?

        \item[(3)] How many $\equiv_{dm}$-equivalence classes of
        $\Sigma^0_4$ quantifier patterns are there?
    \end{enumerate}
\end{problem}

\begin{problem}
    Is there a decision procedure which, given quantifier patterns
    \[
        \overline{P},\overline{Q}
        \in
        \{\exists,\forall,\exists^\infty,\forall^\infty\}^{<\omega},
    \]
    decides whether
    \[
        \langle\overline{P}\rangle
        \leq_m
        \langle\overline{Q}\rangle?
    \]
\end{problem}

\vskip -10pt
\begin{bibentry}{9}

\end{bibentry}
}

\bigskip
\bigskip

\section{Liang Yu}
\abstractbody{

A set $A\subseteq\mathbb{R}^2$ is called a \emph{Kakeya set} if, for every
$\alpha\in\mathbb{R}$, there exists $C_\alpha\in\mathbb{R}$ such that
\[
    y=\alpha x+C_\alpha
    \quad\Longrightarrow\quad
    (x,y)\in A
\]
for every $(x,y)\in\mathbb{R}^2$. Equivalently,
\[
    \bigl\{(x,y)\in\mathbb{R}^2:y=\alpha x+C_\alpha\bigr\}
    \subseteq A
\]
for every $\alpha\in\mathbb{R}$.

\begin{problem}
    Is there a Kakeya set $A\subseteq\mathbb{R}^2$ such that
    \[
        x\oplus y\geq_T\emptyset'
    \]
    for every $(x,y)\in A$?
\end{problem}

A Kakeya set $A$ is called \emph{minimal} if every proper subset
$B\subsetneq A$ is not a Kakeya set.

\begin{problem}
    Does every Borel Kakeya set $A$ contain a minimal Kakeya set
    $B\subseteq A$?
\end{problem}

Yinghe Peng has shown that there is a full fine Kakeya set which contains no minimal Kakeya set.
}

\bigskip
\bigskip

\section{Jialiang He}
These notes concern two related definability questions in infinite combinatorics and Borel model theory.

\begin{definition}
A set $X\subseteq\mathbb R^2$ is called a \emph{two-point set}, if for every line
$L\subseteq\mathbb R^2$,
\[
|X\cap L|=2.
\]
\end{definition}

\begin{theorem}[Larman {\cite{Larman1968}}]
There is no $F_\sigma$ two-point set.
\end{theorem}

Since every $F_\sigma$ set is Borel, Larman's theorem rules out the lowest
nontrivial possibility in the Borel hierarchy, but it does not settle the following problem.

\begin{problem}
Is there a Borel two-point set?
\end{problem}

Miller showed that, assuming G\"odel's axiom of constructibility $V=L$, there is
a $\bm{\Pi}^1_1$, i.e., coanalytic two-point set \cite{Miller1989}.

\begin{definition}
A \emph{standard Borel space} is a measurable space isomorphic to the Borel space of a Polish space. A \emph{Borel field} is a field $(F, +, \times)$ whose universe is a standard Borel space and whose field operations $+,\times$ are Borel.
\end{definition}

\begin{definition}
Let $F$ be a Borel field. A \emph{Borel algebraically closed extension} of $F$ is a Borel field $K$, together with a Borel injective field homomorphism
\[
\iota:F\hookrightarrow K,
\]
such that $K$ is algebraically closed.
\end{definition}

\begin{problem}[Montalb\'an--Nies {\cite[Question 2.8]{MontalbanNies2013}}]
Does every Borel field have a Borel algebraically closed extension?
\end{problem}

\bigskip
\bigskip

\section{Chieu-Minh Tran}
It is about Morley's categoricity in power theorem. 
According to \cite{Mor}, there is the following consequence.
\begin{theorem}
    For a countable theory $T$, if it is $\kappa$ categorial for some $\kappa>\aleph_0$, then it is $\lambda$ categorial for all uncountable $\lambda$.
\end{theorem}
And \cite{CM} replaced $\kappa$ in the above theorem by recursion theoretic complexity (arithmetic degree). Then, the sequel goes to  categoricity without Choice. The first research about this topic date back to \cite{She} shows that we still have the theorem holds if we assume ZF + there is some uncountable well ordered set of reals.
\begin{theorem}{[ZF $+$ there is an uncountable well ordered set
of reals].}

The following conditions on a countable (first order) $T$ are equivalent:
    \begin{enumerate}
        \item $T$ is categorical in some cardinal $\aleph_{\alpha} > \aleph_0$,in $\mathbf V$, of course
        \item $T$ is categorical in every cardinal $\aleph_{\beta} > \aleph_0$,
        in $\mathbf V$, of course
        \item $T$ is, in $\mathbf L[T]$, totally transcendental
        (i.e. $\aleph_0$-stable) with no two cardinal models
        (i.e., for no model $M$ of $T$ and formula
        $\varphi(x,y) \in \mathbf L(\tau_T)$ and
        $\bar a \in {}^{\ell g(\bar y)}M$ do we have
        $\aleph_0 \leq |\varphi(M,\bar a)| < \|M\|$
        and $\|M\|$ is a cardinal, i.e. the set of elements of $M$ is well-orderable hence its power is a cardinal)
        \item If $\mathbf V' \subseteq \mathbf V$ is a transitive class extending
        $\mathbf L$, $T \in \mathbf V'$ and $\mathbf V'$ satisfies ZFC, then the conditions in (3) hold
        \item For some $\mathbf V'$ clause (4) holds.
    \end{enumerate}
\end{theorem}

However, for non-well-ordered cardinals, the situation deteriorates quickly. The sequel is the following.
For such a theory $T$, we want to show it is very close to a stongly minimal $T'$. It is like a vector space over fixed field or algebraically closed field. And since there is no well order, it is hard to take a basis and match them to one another. So there is the question.
\begin{problem}
    For vector space over fixed $\mathbb{K}$, is choose a basis require the full Axiom of Choice?
\end{problem}

\vskip -12pt
\begin{bibentry}{}

\end{bibentry}

\bigskip
\bigskip

\section{Theodore Slaman}
   In 1954, \cite{Mar} shows that for $E\subseteq \mathbb{R}^2$ an analytic set with Hausdorff dimension 1. Then, for almost all $\theta$, the projection of $E$ onto a line with slope $\theta$ has Hausdorff dimension 1. Later in 1969, \cite{Dav} shows that CH could imply there is some $E\subseteq \mathbb{R}^2$ with Hausdorff dimension 1 and for all $\theta$, projection of $E$ onto a line of slope $\theta$ has Hausdorff dimension 0. \par
   \begin{problem}
       Is there a Davies example provably in ZFC?
   \end{problem}
  
  \vskip -10pt
\begin{bibentry}{9}

\end{bibentry}

\bigskip
\bigskip

\section{Hang Zhang}
The following questions are motivated by the problem whether hyper-\linebreak hyperfiniteness implies hyperfiniteness. 

Let $M$ be a   countable model for a sufficiently large
 fragment of ZFC. Suppose that a partial order $P$ is in $M$. Define the space $Gen^P_M$ of $M$-generics for $P$ by 
 
$$Gen^P_M= \{G \subseteq P:  \mbox{ $G$ is an $P$-generic filter over $M$}\}.$$

If we   identify $Gen^P_M$
 as a subspace of $2^P$  with the product topology, then $Gen^P_M$ is $G_\delta$ in  $2^P$ and thus a Polish space. See \cite{smythe} for details. Define $G\sim H $ if and only if $M[G]=M[H]$. Consider the equivalence relation $\sim$ on $Gen^P_M$.
 
\begin{example}
    (Calderoni--Sinapova \cite{CS}) Let $M$ be a countable model for large fragment ZFC.
Suppose $M\models \kappa$ is measurable (witnessed by $\mathcal{U}$). If $M\models \mathbb{P}=\{(s,A): s=\langle s_0,\dots, s_n\rangle,\ n\in\omega,\ A\in\mathcal{U}\}$, with order
$(s,A)\leq (t,B) \iff t\subseteq s \land s\setminus t\subseteq B \land A\subseteq B$.
Then for any generic sequences $\vec{\alpha}, \vec{\beta}$ of $\mathbb{P}$ over $M$, 
\[
M[\vec{\alpha}]=M[\vec{\beta}]\iff \exists N\ \forall n>N\ \ \alpha(n)=\beta(n) \tag{$*$}.
\]
\end{example}

\begin{problem}
    Find more models $M$ and forcing notions $\mathbb{P}\in M$ such that $(*)$ is true.
\end{problem}
\begin{remark}[{\cite[Theorem 3.1]{smythe}}]
If $\mathbb{P}=\mathbb{C}$ (Cohen forcing), then "$x\thicksim y\iff M[x]=M[y]$" is hyper-hyperfinite.
\end{remark}
\begin{problem}[{\cite[Question 5.2]{smythe}}]
 Is $\thicksim$ hyperfinite? 
\end{problem}
\begin{problem}
    Suppose there are a model $M$ and a forcing notion $\mathbb{P}\in M$ such that $(*)$ is true. Can we deduce that $M\models [\exists\  \text{a large cardinal}]$?
\end{problem}
\begin{remark}
    If $$\langle0,\cdots,0,1,0,\cdots,0,1,1,0\cdots\rangle$$
is Cohen, then $$\langle1,\cdots,1,0,1,\cdots,1,0,0,1\cdots\rangle$$ is also Cohen. Notice that these two Cohen reals lead to the same  forcing extension. So Cohen forcing does not satisfy (*) for any $M$.
\end{remark}    

\vskip -10pt
\begin{bibentry}{9}

\end{bibentry}

\bigskip
\bigskip

\section{Daniel Turetsky}
\quad If $\mathbf{a}<\mathbf{b}<\mathbf{c}$ are Turing degrees, and $\mathbf{a}$, $\mathbf{c}$ contain different Martin-L\"of randoms, then must $\mathbf{b}$ contain one? The answer is no.
\begin{problem}
    Does this hold with a larger degree notion, for example, LR-degrees?
\end{problem}

\bigskip
\bigskip

\section{Philip Welch}
    \quad First we recall some basic notions in \cite{Welch}. Let $\beta$ be an ordinal. Suppose $M$ is a nonstandard model of $\mathsf{KP}$, and $\{\zeta_i:i\in\mathbb{N}\}$, $\{s_i:i\in\mathbb{N}\}$ are sequences of $M$-ordinals. If $\left(M,\{\zeta_i:i\in\mathbb{N}\},\{s_i:i\in\mathbb{N}\}\right)$ satisfies the following:
    \begin{enumerate}
        \item $M$ extends $L_\beta$;
        \item for all $i$, $\zeta_i<\zeta_{i+1}<\beta$ and $s_{i+1}<^Ms_i$;
        \item for all $i$, $L_{\zeta_i}\prec_{\Sigma_{m+1}} L_{s_i}^M$;
        \item for all $i$, $L_{s_{i+1}}^M\prec_{\Sigma_{m-1}} L_{s_i}^M$,
    \end{enumerate}
    then we call this triple a strong $\Sigma_{m+1}$-nesting on $\beta$. And we say $\beta$ supports a strong $\Sigma_{m+1}$-nesting if there exist such nonstandard model $M$ of $\mathsf{KP}$ and sequences $\{\zeta_i:i\in\mathbb{N}\},\{s_i:i\in\mathbb{N}\}$ of $M$-ordinals for $L_\beta$.

    \quad We call a function $\Phi:\mathbb{R}\to\mathbb{R}$ an operator ($\mathcal{P}(\mathbb{N})$ and $\mathbb{R}$ are identified). An operator $\Phi$ is said to be 
    \begin{enumerate}
        \item $\bm{\Pi}_m^1$ if there exists a $\bm{\Pi}_m^1$ formula $\phi$ and a parameter $Y$ such that
        $$\forall n\,\forall X\left(n\in\Phi(X)\leftrightarrow\phi(n,X,Y)\right);$$
        \item monotone if $X\subseteq Y$ implies $\Phi(X)\subseteq\Phi(Y)$.
    \end{enumerate}
    The axiom $\bm{\Pi}_{m+1}^1\mathsf{-MI}$ asserts that for every monotone $\bm{\Pi}_{m+1}^1$ operator $\Phi$, there is some prewellordering $W$ such that for all $a\in\textrm{field}(W)=\textrm{dom}(W)\cup\textrm{ran}(W)$, we have 
    $$W_a=\Phi(W_{<a}).$$
    Here
    \begin{align*}
        &W_a=\{b:a<_Wb\},\\
        &W_{<a}=\bigcup\{W_b:b<_Wa\}.
    \end{align*}
    Note that an operator $\Phi$ can be iterated by setting $\Phi^0=\emptyset, \Phi^{<\alpha}=\bigcup_{\gamma<\alpha}\Phi^\gamma$ and $\Phi^\alpha=\Phi(\Phi^{<\alpha})$. Put $\Phi^\infty=\bigcup_\alpha\Phi^\alpha$. We can see that the field of $W$ is $\Phi^\infty$.
\begin{problem}
    If $L_\beta\models\mbox{$\mathit{V}=\mathit{HC}$}$, then
    \begin{align*}
       \left(\mathcal{P}(\omega)\cap L_\beta\models\bm{\Pi}_{m+1}^1\mathsf{-MI}\,\right)\iff\beta\mbox{ supports a strong }\Sigma_{m+1}\textrm{-nesting}.
    \end{align*}
\end{problem}

\vskip -10pt
\begin{bibentry}{9}

\end{bibentry}

\bigskip
\bigskip

\bigskip
\bigskip

\section{Su Gao}
\quad This question is about equivalence relations given by complete first-order theories. The question has been circulating in the community since early 2000s.

Consider a complete first order theory $T$ in a countable language. Let
$\mathrm{Mod}(T)$ be the Polish space of countable models of $T$ with universe $\omega$. Let $\cong\restriction_{\mathrm{Mod}(T)}$ denote the isomorphism relation on $\mathrm{Mod}(T)$. Recall two basic equivalence relations in descriptive set theory:
\begin{enumerate}
    \item $=$ is the identity relation on $\mathbb{R}$;
    \item $=^+$ is the equivalence relation on $\mathbb{R}^\omega$ that $(x_n)_{n\in\omega}=^+(y_n)_{n\in\omega}$ if they represent the same countable subset of reals, i.e., if
    $$ \{x_n\colon n\in\omega\}=\{y_n\colon n\in\omega\}.$$ 
\end{enumerate}

Also note the following fact:
\begin{fact}
    (1) There is a complete theory $T$ such that $=$ is Borel bi-reducible to $\cong\restriction_{\mathrm{Mod}(T)}$. For example, $T$ can be chosen as the theory of the equivalence relation on $\omega$ with infinitely many equivalence classes of any finite size.\\
   (2) There is a complete theory $T'$ such that $=^+$ is Borel bi-reducible to $\cong\restriction_{\mathrm{Mod}(T')}$. For example, $T'$ can be chosen as the theory of infinitely many independent equivalence relations on $\omega$ and every equivalence relation has infinitely many equivalence classes of any finite size.
\end{fact}
So we have the following natural question:
\begin{problem}\label{prob:12.2}
    Does there exist a complete first-order theory $T$ in a countable language such that $=\ <_B\ \cong\restriction_{\mathrm{Mod}(T)}\ <_B\ =^+$?
\end{problem}

\begin{remark}
    Note that $=\ <_B\ =^+$, and there are equivalence relations strictly between them. For example, $E_0,E_\infty,E^\omega_0$. So, they are natural candidates for the complexity of $\cong\restriction_{\mathrm{Mod}(T)}$ for some $T$ if the problem has a positive answer.
\end{remark}

There have been some research (e.g. \cite{Mar07, Koe09, CK10, RS17}) on this problem and more generally on the isomorphism relation $\cong\restriction_{\mathrm{Mod}(T)}$ in the Borel reducibility sense. In particular, Dave Marker showed in \cite{Mar07} that if $S(T)$ (the Stone space of $T$-types) is uncountable, then $T$ is not an example for Problem~\ref{prob:12.2}.

\vskip -12pt 
\begin{bibentry}{9}

\end{bibentry}

\end{document}